\documentclass[reqno]{amsart}

\theoremstyle{plain}
\newtheorem{theorem}{Theorem}
\newtheorem{corollary}{Corollary}
\newtheorem{lemma}{Lemma}
\newtheorem{proposition}{Proposition}

\theoremstyle{remark}
\newtheorem{remark}{Remark}

\theoremstyle{definition}
\newtheorem{definition}{Definition}

\newcommand\resp{resp.\ }
\newcommand\ie{i.e.\ }
\newcommand\cf{cf.\ }

\newcommand\NN{\mathbb N}
\newcommand\ZZ{\mathbb Z}
\newcommand\RR{\mathbb R}

\newcommand\goe{\mathfrak g}
\newcommand\hoe{\mathfrak h}

\newcommand\vep{\varepsilon}
\newcommand\vph{\varphi}

\newcommand\ad{\operatorname{ad}}
\renewcommand\div{\operatorname{div}}
\newcommand\Fl{\operatorname{Fl}}
\newcommand\grad{\operatorname{grad}}
\newcommand\pr{\operatorname{pr}}
\newcommand\tr{\operatorname{tr}}
\newcommand\Tor{\operatorname{Tor}}

\newcommand\ric{\operatorname{ric}}
\newcommand\sym{{\text{\rm sym}}}
\newcommand\tran{{\mathsf{T}}}

\newcommand\vf{\mathfrak X(M)}
\newcommand\vfevol{\mathfrak X_{\text{\rm ex}}(M,\mu)}
\newcommand\vfvol{\mathfrak X(M,\mu)}

\newcommand\vfesym{\mathfrak X_{\text{\rm ex}}(M,\omega)}
\newcommand\vfsym{\mathfrak X(M,\omega)}
\newcommand\vfecsym{\mathfrak X_{\text{\rm ex}}^{\text{\rm c}}(M,\omega)}
\newcommand\vfconfsym{\mathfrak X(M,[\omega])}
\newcommand\Diff{\operatorname{Diff}(M)}
\newcommand\Diffevol{\operatorname{Diff}_{\text{\rm ex}}(M,\mu)}
\newcommand\Diffvol{\operatorname{Diff}(M,\mu)}
\newcommand\Diffesym{\operatorname{Diff}_{\text{\rm ex}}(M,\omega)}
\newcommand\Diffsym{\operatorname{Diff}(M,\omega)}
\begin{document}

\title{Totally geodesic subgroups of diffeomorphisms}
\author{Stefan Haller, Josef Teichmann and Cornelia Vizman}
\address{Stefan Haller, Institute of Mathematics, University of Vienna,
Strudlhofgasse 4, A-1090 Vienna, Austria.}
\email{stefan@mat.univie.ac.at}
\address{Josef Teichmann, Institute of financial and actuarial mathematics, Technical
University of Vienna, Wiedner Hauptstrasse 8--10, A-1040 Vienna,
Austria.} \email{josef.teichmann@fam.tuwien.ac.at}
\address{Cornelia Vizman, West University of Timisoara, Department of Mathematics,
Bd. V.Parvan 4, 1900 Timisoara, Romania.} \email{vizman@math.uvt.ro}
\thanks{This paper was mainly written while the first author stayed
at the Ohio State University. He would also like to thank D.~Burghelea
for raising the question of Proposition~\ref{nofinite} below.}
\thanks{The second author would like to thank Peter W. Michor and Walter
Schachermayer for the perfect working atmosphere, he acknowledges the
support by the research grant FWF Z-36 "Wittgenstein prize" awarded to
Walter Schachermayer.}
\begin{abstract}
We determine the Riemannian manifolds for which the group of exact
volume preserving diffeomorphisms is a totally geodesic subgroup of the
group of volume preserving diffeomorphisms, considering right invariant
$L^2$-metrics. The same is done for the subgroup of Hamiltonian
diffeomorphisms as a subgroup of the group of symplectic
diffeomorphisms in the K\"ahler case. These are special cases of
totally geodesic subgroups of diffeomorphisms with Lie algebras big
enough to detect the vanishing of a symmetric $2$-tensor field.
\end{abstract}
\maketitle

\section{Introduction}

Euler equation for an ideal fluid flow is just the geodesic
equation on the group of volume preserving diffeomorphisms with
right invariant $L^2$-metric, see \cite{A}. In their book Arnold and
Khesin \cite{AK} state, that the subgroup of exact volume
preserving diffeomorphisms is totally geodesic in the group of
volume preserving diffeomorphisms on compact, oriented surfaces.
This means that the group of Hamiltonian diffeomorphisms is
totally geodesic in the group of symplectic diffeomorphisms, which
has a physical interpretation, namely the existence of a single
valued stream function for the velocity field at the initial
moment implies the existence of a single valued stream function at
any other moment of time. We will prove that a closed, oriented surface
having this property either has the first Betti number zero (the Lie algebras
of symplectic and Hamiltonian vector fields coincide) or it is a
flat $2$-torus.

A much more general classification result is actually true:
twisted products of a torus by a Riemannian (\resp K\"ahler)
manifold with vanishing first Betti number are the only Riemannian
(\resp K\"ahler) manifolds, where the exact volume preserving
diffeomorphisms lie totally geodesic in the Lie group of volume
preserving diffeomorphisms (\resp the Hamiltonian diffeomorphisms
in the symplectic diffeomorphisms), see Definition~\ref{twist_prod},
Theorem~\ref{thm_vol} and Theorem~\ref{thm_sym} in section~\ref{results}.

The problem can be formulated in the setting of regular Fr\'echet-Lie
groups: given a regular Fr\'echet-Lie group in the sense of Kriegl-Michor,
see \cite{KM}, and a (bounded, positive definite) scalar product
$g:\goe\times\goe\to\RR$ on the Lie algebra $\goe$, we can define a right
invariant metric on $G$ by
$$
G_x(\xi,\eta)
:=g((T_x\rho^x)^{-1}\xi,(T_x\rho^x)^{-1}\eta)
\quad\text{for $\xi,\eta\in T_xG$,}
$$
where $\rho^x$ denotes the right translation on $G$. The energy
functional of a smooth curve $c:\RR\to G$ is defined by
$$
E(c)
=\int_a^bG_{c(t)}(c'(t),c'(t))dt
=\int_a^bg(\delta^rc(t),\delta^rc(t))dt,
$$
where $\delta^r$ denotes the right logarithmic derivative on the
Lie group $G$.

Assuming $c:[a,b]\to G$ to be a geodesic with respect to the right invariant
(weak) Riemannian metric, variational calculus yields
\begin{align*}
\frac{d}{dt}X_t&=-\ad(X_t)^\tran X_t
\\
X_t&=\delta^rc(t),
\end{align*}
where $\ad(X)^\tran:\goe\to\goe$ denotes the adjoint with respect to the
Hilbert scalar product of $\ad(X)$, see \cite{KM}, which we assume to exist
as bounded linear map $\ad(\cdot)^\tran:\goe\to L(\goe)$. A Lie subgroup
$H\subseteq G$ is totally
geodesic if any geodesic $c$ with $c(a)=e$ and $c'(a)\in\hoe$ stays in $H$.
This is the case if $\ad(X)^\tran X\in\hoe$ for all
$X\in\hoe$. If there is a geodesic in $G$ in any direction
of $\hoe$, then the condition is necessary and sufficient.

The setting for the whole article is the following: Given a
regular Fr\'echet-Lie group $G$ with Lie algebra $\goe$ and
a bounded, positive definite scalar product on $\goe$. We assume, that
$\ad(\cdot)^\tran:\goe\to L(\goe)$ exists and is bounded.
Furthermore we are given a splitting subalgebra $\hoe$, \ie $\hoe$ has an
orthogonal complement $\hoe^\perp$ in $\goe$ with respect to the scalar
product.
We only assume $\goe=\hoe\oplus\hoe^\perp$ as orthogonal direct sum,
in the algebraic sense.
It follows, that $\hoe$ and $\hoe^\perp$ are closed, and that the
orthogonal projections onto $\hoe$ and $\hoe^\perp$ are bounded
with respect to the Fr\'echet space topology.
$\hoe$ is called totally geodesic in $\goe$ if
$\ad(X)^\tran X\in\hoe$ for all
$X\in\hoe$. The following reformulation of the condition
provides the main condition for our work:

\begin{lemma}
In the situation above,
$\hoe$ is totally geodesic in $\goe$ iff
$\langle[X,Y],X\rangle=0$ for all $X\in\hoe$ and $Y\in\hoe^\perp$.
\end{lemma}

In this article we consider the following important examples of
the outlined situation: Let $M$ be a closed, connected and oriented
manifold. The regular Fr\'echet-Lie group $\Diff$ is
modeled on the vector fields $\vf$, a Fr\'echet
space. The Lie algebra is $\vf$ with the negative of the usual Lie
bracket. The symbol $[\cdot,\cdot]$ will denote the usual Lie bracket and
$\ad(X)Y=[X,Y]$. With these conventions, the geodesic
equation on diffeomorphism groups has no minus sign.
The following subgroups are regular
Fr\'echet-Lie subgroups of $\Diff$, see \cite{KM}:

\begin{enumerate}
\item
The group $\Diffvol$ of volume preserving
diffeomorphisms of $(M,\mu)$, where $\mu$ is a volume form on $M$;
its Lie algebra is $\vfvol$ the Lie algebra of
divergence free vector fields.
\item
The group $\Diffevol$ of exact volume preserving
diffeomorphisms of $(M,\mu)$ with Lie algebra
$$
\vfevol
=\{X\in\vf:\text{$i_X\mu$ is an exact differential form}\}.
$$
\item
The group $\Diffsym$ of symplectic diffeomorphisms
of the symplectic manifold $(M,\omega)$ with Lie algebra
$\vfsym$ of symplectic vector fields.
\item
The group $\Diffesym$ of Hamiltonian diffeomorphisms
of $(M,\omega)$ with Lie algebra $\vfesym$ of Hamiltonian vector fields.
\end{enumerate}

Let $(M,g)$ denote a closed connected and orientable Riemannian
manifold with Riemannian metric $g$, $\nabla$ the Levi-Civita
covariant derivative and $\mu$ the canonical volume form on $M$
induced by the metric $g$ and a choice of orientation. By
$\sharp_g:T^*M\to TM$ we denote the geometric lift
$g(\sharp_g\alpha,\cdot)=\alpha$ and by $\flat_g$ its inverse. We
will omit the index $g$ when no confusion is possible. The
Hodge-$*$-operator is given with respect to the volume form
$\mu$ such that $g(\beta,\eta)\mu=\beta\wedge*\eta$ for
$\beta,\eta$ $k$-forms, where $g$ denotes the respective scalar
product on the forms. The exterior derivative is denoted by $d$,
the codifferential by $\delta=(-1)^{n(k+1)+1}*d*$ on $k$-forms.
With this convention $d$ and $\delta$ are adjoint
with respect to scalar product on forms.
Furthermore $\Delta=d\delta+\delta d$.

In the case of $G=\Diff$ the adjoint of $\ad(X)$, with respect to
the induced right invariant $L^2$-structure is given by the expression
$$
\ad(X)^\tran X
=-\nabla_XX-(\div X)X-\frac12\grad(g(X,X)).
$$
We apply here the notions of gradient of a function $\grad f=\flat(df)$
and divergence of a vector field $\div X=-\delta(\flat X)$, \ie
$L_X\mu=\div(X)\mu$. If
$H\subseteq G$ is a subgroup with splitting subalgebra $\hoe\subseteq\goe$,
then for $X\in\hoe$
$$
\ad(X)|^\tran_\hoe X=\pi(\ad(X)|^\tran_\goe X),
$$
where $\pi:\goe\to\hoe$ denotes the orthogonal projection.

In particular we obtain for the volume preserving diffeomorphisms
$\Diffvol$
\begin{align*}
\ad(X)|_{\vfvol}^\tran X &=-\nabla_XX-(\grad p)(X)
\\
\Delta p
&=\div(\nabla_XX).
\end{align*}
Such a function $p$ exists, is unique up to a constant and smooth by
application of the smooth inverse of the Laplacian on its range.
Hence $(\grad p)(X)$ is a well defined smooth vector field.
Remark, that the orthogonal complement to divergence
free vector fields are gradients of some functions, with respect
to the $L^2$-metric, which is easily seen due to the orthogonal
Hodge decomposition
$\vf=\sharp_gd\Omega^0(M)\oplus\sharp_g\ker\delta$, where
$\sharp_g$ denotes the geometric lift. Consequently the geodesic equation
on $\Diffvol$ with right invariant $L^2$-metric is
$$
\frac{d}{dt}X_t
=-\nabla_{X_t}X_t-(\grad p)(X_t),
$$
the Euler equation for an ideal fluid flow, see \cite{A}.

Let $M$ be an closed, connected almost K\"ahler manifold
$(M,g,\omega,J)$, \ie the symplectic form $\omega$, the almost complex
structure $J$ and the Riemannian metric $g$ satisfy the relation
$g(X,Y)=\omega(X,JY)$. Note, that a K\"ahler manifold has a
natural orientation given by $J$. Moreover we have
$\flat_\omega(X)=-(\flat_gX)\circ J=\flat_g(JX)$ and
$\sharp_\omega\vph=-J\sharp_g\vph=\sharp_g(\vph\circ J)$,
especially $\sharp_\omega:T_x^*M\to T_xM$ is an isometry, for $J$ and
$\sharp_g$ are.

For the symplectic diffeomorphisms $\Diffsym$
we obtain the following adjoint
\begin{align*}
\ad(X)|_{\vfsym}^\tran X
&=-\nabla_XX-\frac12\grad(g(X,X))-\sharp_\omega((\delta\alpha)(X))
\\
d((\delta\alpha)(X))
&=-di_{\nabla_XX+\frac12\grad(g(X,X))}\omega
\end{align*}
and the geodesic equation
$$
\frac{d}{dt}X_t
=-\nabla_{X_t}X_t-\frac12\grad(g(X_t,X_t))-\sharp_\omega((\delta\alpha)(X_t)).
$$
Remark that via the symplectic lift the orthogonal Hodge decomposition of
$$
\Omega^1(M)=d\Omega^0(M)\oplus\mathcal{H}^1(M)\oplus\delta\Omega^2(M)
$$
can be carried to the vector fields. Symplectic
vector fields are those with $L_X\omega=d(\flat_\omega X)=0$, \ie
$\flat_\omega X\in\ker d=d\Omega^0(M)\oplus\mathcal{H}^1(M)$. So there is
some (symplectic) harmonic part and some Hamiltonian part, see \cite{KM}.
In the above formula $(\delta\alpha)(X)$ is uniquely determined and smoothly
dependent for $X\in\vfsym$ by the Hodge decomposition. The divergence part
is zero, since symplectic diffeomorphisms are volume preserving.

\section{Statement and proof of the results}\label{results}

In this section we develop the necessary notions and prove the asserted
main results of the article. We shall provide several equivalent
conditions, geometric and analytic ones, in the Riemannian and K\"ahler
case, such that the presented subgroups are totally geodesic.

For a $1$-form $\vph$ we set
$$
(\nabla\vph)^\sym(X,Y)
:=(\nabla_X\vph)(Y)+(\nabla_Y\vph)(X),
$$
the symmetric part of $\nabla\vph$. Note that $\sharp\vph$ is a
Killing vector field, \ie generates a flow of isometries, if and only if
$(\nabla\vph)^\sym=0$. Note also, that
$$
d\vph(X,Y)=(\nabla_X\vph)(Y)-(\nabla_Y\vph)(X),
$$
the skew symmetric part of $\nabla\vph$, and
$\tr(\nabla\vph)^\sym=2\div\sharp\vph=-2\delta\vph$.

\begin{lemma}
Let $(M,g)$ be a closed oriented $n$-dimensional Riemannian manifold and
$\goe\subseteq\vf$ a closed subalgebra, such that
$\ad(X)^\tran:\goe\to\goe$ exists for all $X\in\goe$. Then one has
$$
2\int_Mg(\ad(X)^\tran(X),Y)\mu
=\int_M\big((\nabla\flat Y)^\sym+(\div Y)g\big)(X,X)\mu.
$$
Moreover
$$
\tr\big((\nabla\flat Y)^\sym+(\div Y)g\big)=(n+2)\div Y.
$$
Especially $(\nabla\flat Y)^\sym+(\div Y)g=0$ iff $(\nabla\flat Y)^\sym=0$,
\ie $Y$ is Killing.
\end{lemma}

\begin{proof}
We have:
\begin{align*}
\int g(\ad(X)^\tran(X),Y)\mu
&=\int g(X,[X,Y])\mu
\\
&=\int\big(g(X,\nabla_XY)-g(X,\nabla_YX)\big)\mu
\\
&=\int\big(i_X\nabla_X\flat Y-\tfrac12L_Yg(X,X)\big)\mu
\\
&=\int\tfrac12
(\nabla\flat Y)^\sym(X,X)\mu+\tfrac12g(X,X)L_Y\mu
\\
&=\frac12\int\big((\nabla\flat Y)^\sym+(\div Y)g\big)(X,X)\mu
\end{align*}
The second statement follows from
$\tr\big(\nabla(\flat Y)^\sym\big)=2\div Y$ and $\tr(g)=n$.
\end{proof}

\begin{definition}
We say $\goe\subseteq\vf$ is big enough to detect the
vanishing of a symmetric $2$-tensor field, if a symmetric $2$-tensor field
$T\in\Gamma(S^2T^*M)$ vanishes if
$$
\int_MT(X,Y)\mu=0\quad\text{for all $X,Y\in\goe$.}
$$
\end{definition}

\begin{remark}\label{tg_alt}
Let $\goe\subseteq\vf$ be the Lie algebra of a Lie
group of diffeomorphisms $G$, such that $\ad(\cdot)^\tran:\goe\to L(\goe)$
is bounded.
Suppose $H\subseteq G$ is a Lie subgroup with splitting Lie subalgebra
$\hoe\subseteq\goe$ and assume that $\hoe$ is big enough to detect the
vanishing of a symmetric $2$-tensor. It will then be easy to decide if
$\hoe$ is totally geodesic in $\goe$, for the condition of
being totally geodesic, $\ad(X)^\tran X\in\hoe$
for all $X\in\hoe$, translates in this case to
$(\nabla\flat Y)^\sym=0$ for all $Y\in\hoe^\perp$, \ie
$Y$ is Killing for all $Y\in\hoe^\perp$.

Note also, that if $\hoe$ has an orthogonal
complement in $\vf$, \ie $\hoe\oplus\hoe'=\vf$, then $\hoe$ is a splitting
subalgebra of $\goe$, since $\hoe^\perp=\hoe'\cap\goe$ is the orthogonal
complement of $\hoe$ in $\goe$.
\end{remark}

\begin{lemma}\label{lem_big_sym}
Let $(M,\omega)$ be a symplectic manifold. Then the Lie algebra of compactly
supported Hamiltonian vector fields is big enough to detect the vanishing of
a symmetric $2$-tensor field.
\end{lemma}

\begin{proof}
Suppose $T$ is a symmetric $2$-tensor field, which does not vanish at a point
in $M$. Using Darboux's theorem, and rescaling $\omega$ and $T$ by constants,
we may choose a chart $M\supseteq U\to(-1,1)^{2n}\subseteq\RR^{2n}$, such that
\begin{equation*}
\omega=dx^1\wedge dx^2+\cdots+dx^{2n-1}\wedge dx^{2n}
\end{equation*}
and $T_{22}(x)>0$ for all $x\in(-1,1)^{2n}$, where
$T=\sum T_{ij}dx^i\otimes dx^j$. Now choose a bump
function $b:\RR\to[0,1]$, such that $b(t)=0$ for $|t|\geq\frac12$ and
$b(0)=1$. For $0<\vep\leq1$ we define
\begin{equation*}
\lambda_\vep(x^1,\dotsc,x^{2n})
:=b(\tfrac{x^1}\vep)b(x^2)\cdots b(x^{2n})
\end{equation*}
and $Z_\vep:=\sharp_\omega d\lambda_\vep$. Since the
support of $\lambda_\vep$ is contained in
$(-\vep,\vep)\times(-1,1)^{2n-1}$, $Z_\vep$
extends by zero to a compactly supported Hamiltonian vector field
on $M$. An easy calculation shows
$$
\lim_{\vep\to0}\vep\int_MT(Z_\vep,Z_\vep)\omega^n
=\int_{(-1,1)^{2n}}\big(b'(x^1)b(x^2)\cdots b(x^{2n})\big)^2T_{22}
(0,x^2,\dotsc,x^{2n})\omega^n>0
$$
and hence $\int_MT(Z_\vep,Z_\vep)\omega^n\neq 0$ for
$\vep$ small enough.
\end{proof}

\begin{lemma}\label{lem_big_vol}
Let $(M,\mu)$ be a manifold with volume form, $\dim(M)>1$. Then the Lie
algebra of compactly supported exact divergence free vector fields is big
enough to detect the vanishing of a symmetric $2$-tensor field.
\end{lemma}

\begin{proof}
As in the proof of Lemma~\ref{lem_big_sym}, we choose a chart
$M\supseteq U\to(-1,1)^n\subseteq\RR^n$, such that
$$
\mu=dx^1\wedge\cdots\wedge dx^n
$$
and such that $T_{22}>0$ on $x\in(-1,1)^n$. Take a bump
function $b$ as above and set
$$
\lambda_\vep:=b(\tfrac{x^1}\vep)b(x^2)\cdots b(x^n).
$$
Now define $i_{Z_\vep}\mu:=d(\lambda_\vep dx^3\wedge\cdots\wedge dx^n)$.
Then $Z_\vep$ is a compactly supported exact divergence free vector
field on $M$ and
$$
Z_\vep=b(\tfrac{x^1}\vep)b'(x^2)b(x^3)\cdots b(x^n)
\tfrac{\partial}{\partial x^1}
-\tfrac1\vep b'(\tfrac{x^1}\vep)b(x^2)\cdots b(x^n)
\tfrac{\partial}{\partial x^2}.
$$
Again we get
$$
\lim_{\vep\to 0}\vep\int_MT(Z_\vep,Z_\vep)\mu
=\int_{(-1,1)^n}\big(b'(x^1)b(x^2)\cdots
b(x^n)\big)^2T_{22}(\vep x^1,x^2,\dotsc,x^n)\mu>0,
$$
and hence $\int_MT(Z_\vep,Z_\vep)\mu\neq 0$ for $\vep$ small enough.
\end{proof}

\begin{definition}[Twisted products]\label{twist_prod}
Let $T^k=\RR^k/\Lambda$ be a flat torus, equipped with the metric induced
from the Euclidean metric on $\RR^k$. Suppose $F$ is an oriented
Riemannian manifold and that $\Lambda$ acts on $F$ by orientation
preserving isometries. The total space of the associated fiber bundle
$\RR^k\times_\Lambda F\to T^k$ is an oriented Riemannian
manifold in a natural way. Locally over $T^k$ the metric is the product
metric. We call every manifold obtained in this way a twisted product
of a flat torus and the oriented Riemannian manifold $F$.

If $k$ is even, $F$ K\"ahler and $\Lambda$
acts by isometries preserving the K\"ahler structure then
$\RR^k\times_\Lambda F\to T^k$ is a K\"ahler manifold in a natural way
and we call it a twisted product of a flat torus with the K\"ahler
manifold $F$.
\end{definition}

\begin{theorem}\label{thm_vol}
Let $(M,g)$ be a closed, connected and oriented Riemannian manifold with
volume form $\mu$. Then the following are equivalent:
\begin{enumerate}
\item\label{thm_vol_i}
The group of exact volume preserving diffeomorphisms is a totally geodesic
subgroup in the group of all volume preserving diffeomorphisms.
\item\label{thm_vol_ii}
Every harmonic $1$-form is parallel.
\item\label{thm_vol_iii}
$\ric(\beta_1,\beta_2)=0$ for all harmonic $1$-forms $\beta_1,\beta_2$.
\item\label{thm_vol_iv}
$(M,g)$ is a twisted product of a flat torus and a closed, connected,
oriented Riemannian manifold with vanishing first Betti number.
\item\label{thm_vol_v}
For all $2$-forms $\alpha$ and all harmonic $1$-forms $\beta$ one has
$$
\int_Mg(d\delta\alpha,\delta\alpha\wedge\beta)\mu=0.
$$
\end{enumerate}
\end{theorem}

\begin{proof}[Proof of Theorem~\ref{thm_vol}]
Recall that the orthogonal complement of $\vfevol$ in
$\vfvol$ is $\{\sharp\beta:\text{$\beta$ harmonic
$1$-form}\}$. The equivalence
$(\ref{thm_vol_i})\Leftrightarrow(\ref{thm_vol_ii})$ now follows
immediately from Remark~\ref{tg_alt}, Lemma~\ref{lem_big_vol} and
the fact that for closed $1$-forms $(\nabla\beta)^\sym=0$ is
equivalent to $\nabla\beta=0$.

$(\ref{thm_vol_ii})\Rightarrow(\ref{thm_vol_iii})$ is obvious from the
definition of the curvature
$R_{X,Y}Z=\nabla_X\nabla_YZ-\nabla_Y\nabla_XZ-\nabla_{[X,Y]}Z$ and
$\ric=-\tr_{13}R$.

The integrated Bochner equation on $1$-forms, see for
example \cite{LM}, takes the form
$$
\langle\Delta\alpha,\alpha\rangle
=\|\nabla\alpha\|^2+\int_M\ric(\alpha,\alpha)\mu,
$$
and $(\ref{thm_vol_iii})\Rightarrow(\ref{thm_vol_ii})$ follows.

$(\ref{thm_vol_iv})\Rightarrow(\ref{thm_vol_ii})$: Suppose
$M=\RR^k\times_\Lambda F$. Since $H^1(F;\RR)=0$ it
follows from the Leray-Serre spectral sequence that the projection
$M\to T^k$ induces an isomorphism $H^1(M;\RR)\cong H^1(T^k;\RR)$. So every
harmonic $1$-form comes from $T^k$ and hence is parallel.

$(\ref{thm_vol_ii})\Rightarrow(\ref{thm_vol_iv})$, \cf Theorem~8.5 in
\cite{LM} and \cite{CG}: Suppose $(M,g)$ is a closed,
connected and oriented Riemannian manifold, such that every harmonic
$1$-form is parallel. Choose an orthonormal base
$\{\beta_1,\dotsc,\beta_k\}$ of
harmonic $1$-forms. Since they are parallel they are orthonormal at every
point in $M$. Choose a base point $x_0\in M$, let $U\subseteq\pi_1(M,x_0)$
be the kernel of the Hur\'ewicz-homomorphism
\begin{equation*}
\pi_1(M,x_0)\to H_1(M;\ZZ)\to H_1(M;\ZZ)/\Tor(H_1(M;\ZZ))\cong\ZZ^k,
\end{equation*}
and let $\pi:\tilde M\to M$
be the covering of $M$, which has $U$ as characteristic subgroup. This
is a normal covering, the group of deck transformations is $\ZZ^k$ and
$\pi^*\beta_i=df_i$ for smooth functions $f_i:\tilde M\to\RR$.
Let $z_0$ be a base point in $\tilde M$ sitting above $x_0$ and assume
$f_i(z_0)=0$. Consider the mapping
\begin{equation*}
f:\tilde M\to\RR^k,\quad
f(z)=\bigl(f_1(z),\dotsc,f_k(z)\bigr).
\end{equation*}
Obviously this is a proper, surjective submersion and $F:=f^{-1}(0)$ is a
compact submanifold. Let $X_i:=\sharp_g\pi^*\beta_i$. Then the $X_i$ are
orthonormal at every point and they are all parallel, especially they
commute. Consider
\begin{equation*}
\kappa:F\times\RR^k\to\tilde M,\quad
\kappa(z,t):=\big(\Fl^{X_1}_{t_1}\circ\cdots\circ\Fl^{X_k}_{t_k}\big)(z).
\end{equation*}
Of course we have $f\circ\kappa=\pr_2$, and it follows easily, that $\kappa$
is a diffeomorphism. So $F$ is closed, connected, oriented and
$H^1(F;\RR)=H^1(\tilde M;\RR)=0$. Moreover $\kappa^*g$ is the product metric
of the induced metric on $F$ and the standard metric on $\RR^k$. Every
deck transformation of $\tilde M$ is of the form
\begin{equation*}
F\times\RR^k\to F\times\RR^k,\quad (z,t)\mapsto(\vph_\lambda(z),t+\lambda),
\end{equation*}
where
$\lambda\in\Lambda\subseteq\RR^k$
and $\vph_\lambda$ is an orientation
preserving isometry of $F$. So $M$ is a twisted product, as claimed.

To see $(\ref{thm_vol_i})\Leftrightarrow(\ref{thm_vol_v})$ let
$\sharp\delta\alpha$ be an exact volume preserving vector field,
$\alpha\in\Omega^2(M)$, and let $\beta$ be a harmonic $1$-form. Then
\begin{align*}
\int_Mg(\ad(\sharp\delta\alpha)^{\mathsf{T}}\sharp\delta\alpha,\sharp\beta)\mu
&=\int_Mg(\sharp\delta\alpha,[\sharp\delta\alpha,\sharp\beta])\mu
\\
&=-\int_Mg(\sharp\delta\alpha,\sharp\delta(\delta\alpha\wedge\beta))\mu
\\
&=-\int_Mg(d\delta\alpha,\delta\alpha\wedge\beta)\mu,
\end{align*}
where we used
$$
\delta(\vph_1\wedge\vph_2)-(\delta\vph_1)\wedge\vph_2
+\vph_1\wedge\delta\vph_2
=-\flat[\sharp\vph_1,\sharp\vph_2]
\quad\text{for $\vph_1,\vph_2\in\Omega^1(M)$}
$$
to obtain
$[\sharp\delta\alpha,\sharp\beta]=-\sharp\delta(\delta\alpha\wedge\beta)$
for the second equality.
\end{proof}

\begin{theorem}\label{thm_sym}
Let $(M,g,J,\omega)$ be a closed, connected K\"ahler manifold. Then the
following are equivalent:
\begin{enumerate}
\item\label{thm_sym_i}
The group of Hamiltonian diffeomorphisms is a totally geodesic
subgroup in the group of all symplectic diffeomorphisms.
\item\label{thm_sym_ii}
Every harmonic $1$-form is parallel.
\item\label{thm_sym_iii}
$\ric(\beta_1,\beta_2)=0$ for all harmonic $1$-forms $\beta_1,\beta_2$.
\item\label{thm_sym_iv}
$(M,g,J,\omega)$ is a twisted product of a flat torus and a closed
connected K\"ahler manifold with vanishing first Betti number.
\item\label{thm_sym_v}
For all functions $f$ and all harmonic $1$-forms $\beta$ one has
$$
\int_M(\Delta f)df\wedge\beta\wedge\omega^{n-1}=0.
$$
\end{enumerate}
\end{theorem}

\begin{proof}[Proof of Theorem~\ref{thm_sym}]
Recall that the orthogonal complement of $\vfesym$ in
$\vfsym$ is $\{\sharp_\omega\beta:\text{$\beta$ harmonic
$1$-form}\}$. By Remark~\ref{tg_alt} and Lemma~\ref{lem_big_sym},
(\ref{thm_sym_i}) is equivalent to
$\nabla(\beta\circ J)^\sym=\nabla(\flat_g\sharp_\omega\beta)^\sym=0$
for all harmonic $1$-forms
$\beta$. On a K\"ahler manifold one has $\Delta(\vph\circ
J)=(\Delta\vph)\circ J$ for $1$-forms $\vph$. Particularly
the space of harmonic $1$-forms is $J$-invariant, and so
$(\ref{thm_sym_i})$ is equivalent to $(\nabla\beta)^\sym=0$ and
since harmonic $1$-forms are closed this is equivalent to
(\ref{thm_vol_ii}).

$(\ref{thm_sym_ii})\Leftrightarrow(\ref{thm_sym_iii})$ and
$(\ref{thm_sym_iv})\Leftrightarrow(\ref{thm_sym_ii})$ are as in the proof of
Theorem~\ref{thm_vol}. For $(\ref{thm_sym_iv})\Rightarrow(\ref{thm_sym_ii})$
one needs some extra arguments: One observes,
that the span of the $X_i$ constructed in the proof of
Theorem~\ref{thm_vol},
is $J$-invariant and so is its orthogonal complement. Hence $F$ is a
complex submanifold and therefore a K\"ahler submanifold. Moreover
the complex structure is, locally over $T^k$, the product structure
and so is the symplectic structure as well.

$(\ref{thm_sym_i})\Leftrightarrow(\ref{thm_sym_v})$ follows from
the following computation for a function $f$ and a closed $1$-form $\beta$:
\begin{align*}
\int_Mg(\ad(\sharp_\omega df)^{\mathsf{T}}\sharp_\omega df,\sharp_\omega\beta)\omega^n
&=\int_Mg(\sharp_\omega df,[\sharp_{\omega}df,\sharp_\omega\beta])\omega^n
\\
&=-\int_Mg(\sharp_\omega df,\sharp_\omega d(L_{\sharp_\omega\beta}f))\omega^n
\\
&=-\int_Mg(df,d(L_{\sharp_\omega\beta}f))\omega^n
\\
&=-\int_M(\Delta f)(L_{\sharp_\omega\beta}f)\omega^n
\\
&=-n\int_M(\Delta f)df\wedge\beta\wedge\omega^{n-1}
\end{align*}
For the second equality we used
$[\sharp_\omega\vph_1,\sharp_\omega\vph_2]
=-\sharp_\omega(L_{\sharp_\omega\vph_2}\vph_1)$
for closed $1$-forms $\vph_1,\vph_2$, a relation derived from
$i_{[X,Y]}=L_Xi_Y-i_YL_X$.
\end{proof}

\begin{remark}
The fact that $M$ is K\"ahler was only used to show, that the space of
harmonic $1$-forms is invariant under $J$. In the almost K\"ahler case
the arguments in the proof of Theorem~\ref{thm_sym} show, that
following are equivalent:
\begin{enumerate}
\item
The group of Hamiltonian diffeomorphisms is a totally geodesic subgroup
in the group of all symplectic diffeomorphisms.
\item
$\sharp_\omega\beta=\sharp_g(\beta\circ J)=-J\sharp_g\beta$ is Killing
for every harmonic $1$-form $\beta$.
\item
For all functions $f$ and all harmonic $1$-forms $\beta$ one has
$$
\int_M(\Delta f)df\wedge\beta\wedge\omega^{n-1}=0.
$$
\end{enumerate}
\end{remark}

\begin{remark}
The computation in the proof of Theorem~\ref{thm_sym} shows, that for a
function $f$, $\ad(\sharp_\omega df)^\tran\sharp_\omega df=0$ if and
only if $$ \int_M(\Delta f)df\wedge\beta\wedge\omega^{n-1}=0
\quad\text{for all closed $1$-forms $\beta$,} $$ even on almost
K\"ahler manifolds. If $f$ is a `generalized eigenvector' of the
Laplacian, \ie $\Delta f=h\circ f$ for some smooth function $h\in
C^{\infty}(\RR,\RR)$, then $$ \int_M(\Delta
f)df\wedge\beta\wedge\omega^{n-1} =\int_M(h\circ
f)df\wedge\beta\wedge\omega^{n-1} =\int_{M}d\bigl((H\circ
f)\beta\wedge\omega^{n-1}\bigr)=0 $$ with $H$ an integral of $h$,
consequently the condition is satisfied. So the geodesic is given by an
exponential. These are examples of a more general method how to solve
the geodesic equation: In the general setting any finite dimensional
submanifold $S\subset\goe$ such that $\ad(X)^\tran X\in T_XS$
for $X\in S$ admits the calculation of flowlines in the manifold $S$.
In the above case the submanifold $S$ is given by a single point.
\end{remark}

\section{Final Remarks}

Since Killing vector fields are divergence free, Remark~\ref{tg_alt}
immediately implies

\begin{corollary}
Let $M$ be a closed, connected and oriented Riemannian manifold. Then
there does not exist a closed Lie subalgebra
$\vfvol\subset\goe\subseteq\vf$, such that
$\ad(\cdot)^\tran:\goe\to L(\goe)$ is bounded and such that
$\vfvol$ is totally geodesic in $\goe$.
\end{corollary}

On a closed Riemannian manifold the Lie algebra of Killing vector fields is
finite dimensional. So Remark~\ref{tg_alt} also implies

\begin{corollary}
Let $(M,g,J,\omega)$ be a closed, connected almost K\"ahler manifold. Then
the symplectic diffeomorphisms are not totally geodesic in the group of
volume preserving diffeomorphisms, provided $\dim(M)>2$.
\end{corollary}

\begin{remark}
Let $K$ be a compact Lie group acting by isometries on the closed connected
orientable manifold $(M,g)$. In \cite{V} it is shown that the group of
$K$-equivariant diffeomorphisms is a totally geodesic subgroup of
$\Diff$. Its Lie algebra, $K$-invariant vector fields on $M$,
is split, a complement is $\{X\in\vf:\int_Kk^*Xdk=0\}$, infinite
dimensional. This does not contradict the arguments above, for
the Lie algebra of $K$-invariant vector fields on $M$ does not detect the
vanishing of a symmetric $2$-tensor field.
\end{remark}


Let $\vfecsym$ denote the compactly supported Hamiltonian vector fields.

\begin{lemma}
Let $(M,\omega)$ be a connected symplectic manifold and let
$\alpha\in\Omega^2(M)$. If
$$
L_X\alpha=0,\quad\text{for all $X\in\vfecsym$}
$$
then there exists $\lambda\in\RR$ such that $\alpha=\lambda\omega$.
\end{lemma}

\begin{proof}
Choose a Darboux chart centered at $z\in M$, such that
\begin{equation*}
\omega=dx^1\wedge dy^1+\cdots+dx^n\wedge dy^n
\end{equation*}
and write
\begin{equation*}
\alpha=\sum_{i<j}a_{ij}dx^i\wedge dx^j+\sum_{i<j}b_{ij}dy^i\wedge dy^j+
\sum_{i,j}c_{ij}dx^i\wedge dy^j.
\end{equation*}
Let $h$ be a compactly supported function on $M$, such that
$h=x^i$ \resp $h=y^i$ locally around $z$. Then the condition
$L_{\sharp_\omega dh}\alpha=0$ shows that $a_{ij}$, $b_{ij}$ and
$c_{ij}$ are all constant locally around $z$. Using $h=(x^i)^2$
one sees, that $b_{ij}=0$ and $c_{ij}=0$ for $i\neq j$. Using
$h=(y^i)^2$ yields
$a_{ij}=0$. Finally, using $h=x^ix^j$ shows $c_{ii}=c_{jj}$. So
$\alpha=\lambda\omega$ locally around $z$, for some constant
$\lambda\in\RR$. Since $M$ is connected this is true globally.
\end{proof}

We denote by
$$
\vfconfsym:=\{X\in\vf:\exists\lambda\in\RR:L_X\omega=\lambda\omega\}.
$$
Notice that for closed $M$ we have $\vfconfsym=\vfsym$. Moreover if
$L_X\omega=f\omega$ for some function $f$ and if $\dim(M)>2$ then $f$ is
constant, due to the non-degeneracy of $\omega$.

\begin{lemma}\label{lastlem}
Let $(M,\omega)$ be a symplectic manifold and let $Z\in\vf$. If
$$
[Z,X]\in\vfconfsym,\quad\text{for all $X\in\vfecsym$}
$$
then $Z\in\vfconfsym$.
\end{lemma}

\begin{proof}
Set $\alpha:=L_Z\omega\in\Omega^2(M)$. Then for every $X\in\vfecsym$ we
have
\begin{equation*}
L_X\alpha=L_XL_Z\omega=L_{[X,Z]}\omega=\lambda\omega=0.
\end{equation*}
Here $\lambda$ has to vanish, since $[X,Z]$ has compact support. So by
the previous lemma there exists $\tilde\lambda\in\RR$, such that
$L_Z\omega=\alpha=\tilde\lambda\omega$, \ie
$Z\in\vfconfsym$.
\end{proof}

\begin{proposition}\label{nofinite}
Let $(M,\omega)$ be a symplectic manifold. Then there does not exists a Lie
subalgebra $\vfconfsym\subset\goe\subseteq\vf$, such that
$\vfconfsym$ has finite codimension in $\goe$.
\end{proposition}

\begin{proof}
Suppose $\goe$ is bigger than $\vfconfsym$. Then there exists $Z\in\goe$
and an open subset $U\subseteq M$, such that
$Z|_V\notin\mathfrak X(V,[\omega])$, for every
all open $V\subseteq U$. For any $k\in\NN$ we choose disjoint subsets
$V_1,\dotsc,V_k\subseteq U$. Since
$Z|_{V_i}\notin\mathfrak X(V_i,[\omega])$ Lemma~\ref{lastlem}
yields $X_i\in\mathfrak X_{\text{\rm ex}}^{\text{\rm c}}(V_i,\omega)$, such
that $Y_i:=[Z,X_i]\notin\mathfrak X(V_i,[\omega])$. But $Y_i\in\goe$ and
obviously $\{Y_1,\dotsc,Y_k\}$ are linearly independent in
$\goe/\vfconfsym$. Hence the codimension of
$\vfconfsym$ in $\goe$ is at least $k$. Since $k$ was
arbitrary we are done.
\end{proof}

\begin{corollary}
Let $(M,g,J,\omega)$ be a closed, connected almost K\"ahler manifold.
Then there does not exist a closed Lie subalgebra
$\vfsym\subset\goe\subseteq\vf$, such that $\ad(\cdot)^\tran:\goe\to
L(\goe)$ is bounded and such that $\vfsym$ is totally geodesic in
$\goe$.
\end{corollary}

\begin{proof}
Suppose conversely such a $\goe$ exists. By Remark~\ref{tg_alt} and
Lemma~\ref{lem_big_sym}, $\vfsym$ has
an orthogonal complement in $\goe$, consisting of Killing vector fields. So
this complement has to be finite dimensional, but this contradicts
Proposition~\ref{nofinite}.
\end{proof}

For a manifold with volume form $(M,\mu)$ we let
$$
\mathfrak X(M,[\mu]):=\{X\in\mathfrak X(M):\exists\lambda\in:L_X\mu
=\lambda\mu\}.
$$
Notice, that for closed $M$ one has $\vfvol=\mathfrak X(M,[\mu])$.
Similarly, although it does not yield anything new for our
totally geodesic subgroups, one shows

\begin{proposition}
Let $(M,\mu)$ be a manifold with volume form and $\dim(M)>1$. Then there does
not exist a Lie subalgebra $\mathfrak X(M,[\mu])\subset\goe\subseteq\vf$,
such that $\mathfrak X(M,[\mu])$ has finite codimension in $\goe$.
\end{proposition}

\end{document}